\documentclass{amsart}
\usepackage{amssymb, amsmath, amscd, latexsym, mathrsfs, appendix}
\usepackage{stmaryrd}
\usepackage{graphics}
\usepackage{amsmath}
\usepackage{amsxtra}
\usepackage{amsfonts}
\usepackage{appendix}
\usepackage[all]{xy}
\usepackage{tikz}
\usetikzlibrary{matrix,arrows}

\usepackage{color}

\numberwithin{equation}{section}
\newtheorem{thm}{Theorem}[section]
\newtheorem{cor}[thm]{Corollary}
\newtheorem{lem}[thm]{Lemma}

\newtheorem{prop}[thm]{Proposition}

\theoremstyle{definition}
\newtheorem{defn}[thm]{Definition}

\newtheorem{rem}[thm]{Remark}

\newcommand{\lra}{\longrightarrow}

\newcommand{\co}{\colon\!}

\newcommand{\holim}{\textup{holim}}

\newcommand{\hofiber}{\textup{hofiber}}

\newcommand{\ob}{\textup{ob}}
\newcommand{\map}{\textup{map}}
\newcommand{\mapcov}{\textup{mapcov}}

\newcommand{\rmap}{\mathbb R\textup{map}}

\newcommand{\latch}{\textup{latch}}

\newcommand{\topcat}{\mathsf{Top}}

\newcommand{\emb}{\textup{emb}}
\newcommand{\config}{\mathsf{con}} 

\newcommand{\loc}{\textup{loc}}
\newcommand{\imm}{\textup{imm}}

\newcommand{\fin}{\mathsf{Fin}}
\newcommand{\epifin}{\mathsf{biFin}}
\newcommand{\trifin}{\mathsf{triFin}}
\newcommand{\epifinone}{\mathsf{A}}

\newcommand{\inj}{\mathsf{Inj}}

\newcommand{\sC}{\mathcal C}

\newcommand{\sE}{\mathcal E}

\newcommand{\sP}{\mathcal P}

\newcommand{\op}{\textup{op}}

\newcommand{\RR}{\mathbb R}

\newcommand{\colimsub}[1]{\begin{array}[t]{cc} \textup{colim} \\
[-1.7mm] \scriptstyle{#1} \end{array}}

\newcommand{\holimsub}[1]{\begin{array}[t]{cc} \textup{holim} \\ [-1mm]
\scriptstyle{#1} \end{array}}

\newcommand{\hocolimsub}[1]{\begin{array}[t]{cc} \textup{hocolim} \\
[-1.2mm] \scriptstyle{#1} \end{array}}

\newcommand{\uli}{\underline}

\newcommand{\downmap}[2]{\vcenter{\halign{
&\hfil$##$\hfil\cr#1\cr\downarrow\cr#2\cr}}
}

\begin{document}
\pdfoutput=1   

\title{The configuration category of a covering space}
\author{Pedro Boavida de Brito and Michael S. Weiss}%

\address{Dept. of Mathematics, Instituto Superior Tecnico, Univ.~of Lisbon, Av.~Rovisco Pais, Lisboa, Portugal}%
\email{pedrobbrito@tecnico.ulisboa.pt}

\address{Math. Institut, Universit\"at M\"{u}nster, 48149 M\"{u}nster, Einsteinstrasse 62, Germany}%
 \email{m.weiss@uni-muenster.de}

\thanks{The project was funded by the Deutsche Forschungsgemeinschaft (DFG, German Research Foundation) – Project-ID 427320536 – SFB 1442,
as well as under Germany’s Excellence Strategy EXC 2044 390685587, Mathematics Münster: Dynamics–Geometry–Structure.
P.B. was supported by FCT 2021.01497.CEECIND and grant SFRH/BPD/99841/2014.}

\subjclass[2000]{57R40, 55U40, 55P48}
\begin{abstract} We investigate the relationship between the configuration category
of a manifold and the configuration category of a covering space of that manifold.
\end{abstract}
\maketitle


\section{Introduction} \label{sec-intro}
The configuration category $\config(M)$ of a manifold $M$ organizes the ordered configuration spaces of $M$
for the various finite cardinalities into a bigger homotopical entity. It is a Segal space over the nerve of $\fin$,
where $\fin$ is a standard skeleton of the category of finite sets. See \cite{BoavidaWeissLong} for
a definition, in fact many equivalent definitions. The idea goes back to Andrade's thesis \cite{Andrade}, if not further; closely related abstractions emerged in
connection with the development of factorization homology.

In \cite{BoavidaWeissLong} we emphasized the usefulness of configuration categories in manifold calculus, specifically
in manifold calculus applied to the homotopy theory of spaces of smooth embeddings.
An embedding of manifolds, $L\to M$, determines a map of Segal spaces $\config(L)\to \config(M)$
over the nerve of $\fin$. The resulting map from $\emb(L,M)$ (space of smooth embeddings, where applicable)
or from $\emb^t(L,M)$ (space of topological embeddings) to the space of derived maps
\[  \rmap_\fin(\config(L),\config(M))  \]
is a good starting point for a homotopical analysis of these embedding spaces.

Once it is accepted that configuration categories of manifolds are a good thing, questions of a certain type
arise, namely: does this or that construction with manifolds have an analogue for configuration categories.
In \cite{BoavidaWeissconftensor} we investigated the relationship between the configuration
category of a product of manifolds and the configuration categories of the factors. Here we begin with a
covering space $\pi\co E\to M$ and describe the relationship between $\config(E)$ and $\config(M)$.
These two themes are closely related since a covering space is a special case of a bundle, i.e.,
of a twisted product. As a consequence this paper has much in common with \cite{BoavidaWeissconftensor}.

\subsection{Statement of results}
We begin with topological manifolds $L$ and $M$ and covering spaces
\[ \pi_L\co E_L\to L~,\qquad \pi_M\co E_M\to M\,. \]
Let $\mapcov(\pi_L,\pi_M)$ be the space of pairs $(f,g)$ where $f\co L\to M$ is any (continuous) map,
$g\co E_L\to E_M$ is a map such that $\pi_M g=f\pi_L$~, and for every $z\in L$ the map
\[  \pi_L^{-1}(z) \lra \pi_M^{-1}(f(z))  \]
obtained by restricting $g$ is injective. There is a forgetful map
\begin{equation} \label{eqn-qproj}
q\co \mapcov(\pi_L,\pi_M) \lra \map(L,M)
\end{equation}
which is a covering space.
By making a base change (pullback) along $q$, we obtain
\[  q^*\emb^t(L,M)\lra \emb^t(L,M) \]
which is again a covering space. Elements of $q^*\emb^t(L,M)$ can be interpreted as pairs $(f,g)$
where $f\co L\to M$ is a (topological) embedding a.k.a.~injective continuous map, $g\co E_L\to E_M$ is also a topological embedding and $\pi_M g=f\pi_L$\,.
In particular there is a forgetful map
\begin{equation} \label{eqn-watchout} q^*\emb^t(L,M) \lra \emb^t(E_L,E_M)  \end{equation}
given by $(f,g)\mapsto g$ in the above interpretation.

We aim for a similar picture involving spaces of derived maps between configuration
categories. There is a forgetful map
\[ \rmap_\fin(\config(L),\config(M)) \lra \rmap(L,M)=\map(L,M) \]
obtained by evaluating derived maps of configuration categories on the configurations of cardinality 1.
Using this map, we can make a base change to obtain
\[ q^*\rmap_\fin(\config(L),\config(M)) \lra \rmap_\fin(\config(L),\config(M))  \]
which is again a covering space. (Our mapping space $\rmap_\fin(\config(L),\config(M))$ is actually
a simplicial set. Therefore we allow ourselves to think of $q$ as a map of simplicial sets, by applying the
functor ``singular simplicial set''.) Our main result is
\begin{thm}\label{thm-main} There is a dashed arrow
\[
\xymatrix@C=30pt{
q^*\emb^t(L,M) \ar[d] \ar[r]^-{(\ref{eqn-watchout})} & \emb^t(E_L,E_M) \ar[d] \\
q^*\rmap_{\fin}(\config(L),\config(M)) \ar@{..>}[r] & \rmap_{\fin}(\config(E_L),\config(E_M))
}
\]
making the square homotopy commute.
\end{thm}
The four terms are viewed as simplicial sets.
The vertical arrows in the square are (supposed to be) obvious or nearly so. For the left-hand column,
use the observation that the standard map from $\emb^t(L,M)$ to $\rmap_{\fin}(\config(L),\config(M))$ is a map
over $\map(L,M)$, and $q^*$ is a functor from spaces over $\map(L,M)$ to spaces.
The construction of the dashed arrow is not direct, but it has good functorial properties.
(It is a natural zig-zag $\rightarrow\cdot \leftarrow$ where the second arrow is a weak equivalence.)

\medskip
Here is a summary of the paper. From a covering space $\pi\co E \to M$, we produce a map
$\config(\pi) \to N\fin \times N\fin$ of Segal spaces.
This construction comes with forgetful maps
\[
\config(E) \gets \config(\pi) \to \config(M)
\]
over the projections $N\fin \gets N\fin \times N\fin \to N\fin$. The definition and its basic properties are the theme of
sections \ref{sec-boxes} and \ref{sec-pathlifting}. In section \ref{sec-break}, we show that the left-hand map is a weak equivalence
over $N\fin$, for a certain meaning of weak equivalence. In section~\ref{sec-lifting}, we prove theorem \ref{thm-main} by a higher form of diagram chase and in
section~\ref{sec-loc} we discuss localization, closely related to homotopy sheafification. In section~\ref{sec-cohere} we make a few remarks
on towers of covering spaces.

\section{Lifting configurations} \label{sec-boxes}
As in \cite{BoavidaWeissLong} and \cite{BoavidaWeissconftensor} we write $\fin$ for a standard skeleton of the category
of all finite sets. It has objects $\uli k=\{1,2,\dots,k\}$ for $k\ge 0$. All maps between these sets qualify as morphisms.
For better or worse, we tend to make no distinction between the category $\fin$ and its nerve, a simplicial set.

As in \cite{BoavidaWeissconftensor} we shall say that a surjective map
$f\co \uli k\to \uli \ell$ is \emph{selfic} if the injective map from $\uli\ell$ to $\uli k$ defined by
$j\mapsto \min(f^{-1}(j))$ is monotone. The set of selfic maps from $\uli k$ to $\uli\ell$ is in canonical
bijection with the set of equivalence relations on $\uli k$ having exactly $\ell$ equivalence classes.

\begin{defn} \label{defn-epifin0}
Let $\epifin$ be the following category. An object of $\epifin$ is a selfic (surjective) map
$\uli k\to \uli \ell$ in $\fin$. A morphism in $\epifin$ is a commutative square
\[
\xymatrix@R=12pt@M=6pt{
\uli k \ar[d] \ar[r] & \uli \ell \ar[d] \\
\uli m \ar[r] & \uli n
}
\]
in $\fin$ where both horizontal arrows are selfic; the top row is the source and the bottom row is the target.
There is one more condition: the induced map from $\uli k$ to the pullback $\uli m\times_{\uli n} \uli\ell$
is required to be injective.
\end{defn}

Let $M$ be a topological manifold and let $\pi\co E\to M$ be a covering space, so that $E$ is also a manifold.
The configuration category $\config(M)$ is a simplicial space, but we can view it as the nerve of a certain
topological category whose space of objects is the disjoint union of the spaces $\emb(\uli k\,,M)$ for $k\ge 0$.
In \cite{BoavidaWeissLong} this incarnation of $\config(M)$ is called the \emph{particle model}. --- The following is a hybrid between
$\config(E)$ and $\config(M)$, related to both by forgetful maps.

\begin{defn} \label{defn-epifin1} Let $\config(\pi)$ be the (nerve of the) following topological category. An object of $\config(\pi)$
consists of an object $f\co \uli k\to E$ in $\config(E)$, an object $g\co \uli\ell\to M$
in $\config(M)$ and an object $p\co \uli k\to \uli\ell$ in $\epifin$
such that $\pi f=gp$.
\[
\xymatrix@R=16pt{
\uli k \ar[d]_-f \ar[r]^-p & \uli\ell \ar[d]^-g \\
E \ar[r]^-\pi & M
}
\]
(Note that $p$ and $g$ are both determined by $f$.) A morphism of $\config(\pi)$
with the above source $(f,g,p)$ and target $(f',g',p')$
consists of a morphism in $\config(E)$ from $f$ to $f'$, a morphism in $\config(M)$
from $g$ to $g'$ and a morphism in $\epifin$ from $p$ to $r$ which
are similarly related. In particular the morphism $f\to f'$ involves a homotopy
$(h_t)_{t\in[0,a]}$ and the morphism $g\to g'$ involves a homotopy $(h'_t)_{t\in[0,a]}$,
with the same $a$; we must have $\pi h_t=h'_tq$ for all $t\in [0,a]$.
Consequently a morphism in $\config(\pi)$ is determined by its source, its target and the underlying morphism in $\config(M)$.

The object space of $\config(\pi)$ is topologized as a subspace of the object space of the product
$\config(E)\times \config(M)$; similarly, the morphism space of $\config(\pi)$ can be identified with
a subspace of the morphism space of $\config(E)\times \config(M)$.

(Another way to describe the
object space of $\config(\pi)$ is as follows. The object space of $\config(E)$ is a stratified space. The
strata correspond to pairs $(k,\ell)$ of non-negative integers where $k\ge \ell$,
specifying the cardinality of configurations in $E$ and their
images in $M$. The object space of $\config(\pi)$ is the coproduct of the strata of the object space of $\config(E)$.)

There are evident forgetful maps as follows:
\begin{itemize}
\item[-] a reference map from $\config(\pi)$ to $\epifin$;
\item[-] a forgetful map from $\config(\pi)$ to $\config(E)$ which covers the
functor \emph{source} from $\epifin$ to $\fin$;
\item[-] a forgetful map from $\config(\pi)$ to $\config(M)$ which
covers the functor \emph{target}
from $\epifin$ to $\fin$.
\end{itemize}
We note that $\config(\pi)$ is a Segal space, since it is the nerve of a topological category
with some good properties (e.g., the map \emph{source} from the space of morphisms to the space of objects is a
fibration). The reference map from $\config(\pi)$ to $\epifin$ makes $\config(\pi)$ a fiberwise complete
Segal space over $\epifin$.
\end{defn}

\section{Unique path lifting from a categorical viewpoint}\label{sec-pathlifting}
Here we use elementary path lifting arguments in order to learn something
about the forgetful map $\config(\pi)\to \config(M)$.

\begin{defn} \label{defn-finoverspace} Let $\fin(M)$ be the (nerve of the) following topological category.
An object is an
object $\uli k$ of $\fin$ together with a map $\uli k \to M$ (which is \emph{not} required to be injective).
The space of objects is therefore
\[  \coprod_{k\ge 0} M^k \,. \]
A morphism in $\fin(M)$ has a source $f\co \uli k\to M$ and a target $g\co\uli\ell\to M$;
apart from that it consists of a map $u\co\uli k\to \uli\ell$ and a (Moore) homotopy from
$f$ to $gu$. Therefore $\fin(M)$ contains $\config(M)$. There is a forgetful (continuous)
functor from $\fin(M)$ to $\fin$. Also, $\fin$ can be identified with $\fin(*)$.

Next, let $\fin(\pi)$ be the (nerve of the) following category. An object is an object of $\epifin$,
say $p\co \uli k\to \uli\ell$\,, together with a map from $p$ to $\pi$, that is to say,
a commutative diagram
\[
\xymatrix@R=14pt{
\uli k \ar[d] \ar[r]^-p & \uli\ell \ar[d] \\
E \ar[r]^-\pi &  M
}
\]
This looks exactly like the diagram in definition~\ref{defn-epifin1}, but the vertical arrows are not required to be
injective. Instead we impose a weaker condition: the induced map from $\uli k$ to $E\times_M\uli\ell$ is injective.
A morphism in $\fin(\pi)$ consists of two objects in $\fin(\pi)$ (source and target, denoted $x$ and $y$ for short),
a morphism $u$ in $\epifin$ relating the corresponding objects in $\epifin$, and a Moore path in the object space of $\fin(\pi)$
from $x$ to $u^*y=y\circ u$. Therefore $\fin(\pi)$ contains $\config(\pi)$.
\end{defn}

\begin{prop} \label{prop-configcov} The commutative square of simplicial spaces
\[
\xymatrix{ \config(\pi) \ar[r]^-{\textup{incl.}} \ar[d]_-{\textup{functor \emph{target}}} &
\fin(\pi)\ar[d]^-{\textup{functor \emph{target}}} \\
\config(M)  \ar[r]^-{\textup{incl.}} & \fin(M)
}
\]
is both a homotopy pullback square and a strict pullback square.
\end{prop}

\proof The four maps have been defined, the square is commutative and it is a strict pullback square by inspection. The right-hand
vertical arrow is, in each degree, a covering map, therefore a fibration. \qed

\medskip
Recall that for a Segal
space $X$ and an element $y\in X_0$~, the comma construction $(X\downarrow y)$ is a Segal space defined as follows:
\[  (X\downarrow y)_r:= \hofiber_y[d_1d_2\cdots d_{r+1}\co X_{r+1} \to X_0]. \]
In words: the degree $r$ part of $(X\downarrow y)$ is the homotopy fiber over $y$ of the \emph{ultimate target} operator
from $X_{r+1}$ to $X_0$. For example, if $X=N\sC$ for a (discrete) small category $\sC$, then $y$ is an object of $\sC$
and $(X\downarrow y)$ is exactly the nerve of the comma construction $(\sC\downarrow y)$.

An important feature of configuration categories is the local nature of their comma constructions.
Let $M$ be a manifold, $U$ an open subset of $M$ and $y$ an ordered configuration in $U$,
in short $y\in \config(U)_0$\,. Then the inclusion
\[  (\config(U)\downarrow y) \lra (\config(M)\downarrow y)  \]
is a (degreewise) weak equivalence. This is \cite[cor.~3.6]{BoavidaWeissLong}. Moreover, if $U$ happens to be a
tubular neighborhood in $M$ of the configuration $y$, then $y$ is a weakly terminal object in $\config(U)$ and
so the forgetful map
\[  (\config(U)\downarrow y) \lra \config(U) \]
is also a weak equivalence. Therefore, if $U$ and $y$ are in this special relationship,
we arrive at a zigzag of weak equivalences connecting  $(\config(M)\downarrow y)$ to $\config(U)$.

\begin{cor} \label{cor-configcov} Let $U$ be an open subset of $M$. Let $y$ be an element of $\config(\pi|_U)_0$ and let
$z\in \config(U)_0$ be the value of $y$ under the forgetful map $\config(\pi|_U)\to \config(U)$.
The inclusion of comma constructions
\[  (\config(\pi|_U) \downarrow y) \lra (\config(\pi)\downarrow y) \]
is a weak equivalence. If $U$ happens to be a tubular neighborhood in $M$ of the ordered configuration $z$,
then the forgetful map
\[  (\config(\pi|_U) \downarrow y) \lra \config(\pi|_U)  \]
gives a degreewise weak equivalence of $(\config(\pi|_U) \downarrow y)$ with a full Segal subspace of $\config(\pi|_U)$
determined by a union of selected path components of $\config(\pi|_U)_0$.
\end{cor}

\proof We apply proposition~\ref{prop-configcov}, both with $\pi\co E\to M$ and
with $\pi|_U\co E|_U\to U$. The object $y$ of the Segal space in the top left-hand term
of the commutative square of the proposition determines elements in degree 0 for the other three terms
which we still denote by $y$ for the top row, and by $z$ for the bottom row.
Taking the comma construction for each of these give another homotopy pullback square. Therefore, for the first
assertion, it suffices to show that the inclusions
\[  (\fin(\pi|_U)\downarrow y) \to (\fin(\pi)\downarrow y), \]
\[   (\fin(U)\downarrow z) \lra (\fin(M)\downarrow z) \]
are (degreewise) weak equivalences. But this is clear.

The proof of the second assertion is similar. We are assuming that $U$ is a tubular neighborhood of $z$.
We need to show that
\begin{itemize}
\item[(i)] $y$ is a weakly terminal object for a full Segal subspace of $\fin(\pi|_U)$ determined
by a selection of path components of $(\fin(\pi|_U))_0$;
\item[(ii)] $z$ is a weakly terminal object for $\fin(U)$.
\end{itemize}
Again, this is clear. (The path components to be selected in (i) are those containing objects whose
value in $\fin(E)_0$ meets only path components of $E|_U$ which are also met by the value of $y$ in $\fin(E)_0$.)
\qed

\smallskip
Next we want to say a few things about covering spaces in relation to derived maps between simplicial
spaces. Here we have to make a choice of model category setting. Write $\topcat$ for the
category of compactly generated weak Hausdorff spaces. It will
be most convenient for us to equip this with the model category structure where the fibrations are the
classical Hurewicz fibrations and the weak equivalences are the weak homotopy equivalences.
See \cite{Cole2006}, where useful characterizations of the cofibrations can be found.
This is a simplicial model structure. For the category of simplicial spaces
(simplicial objects in $\topcat$) we choose the Reedy model structure. Therefore a map $f\co X\to Y$
of simplicial spaces is a \emph{weak equivalence} if each $f_n\co X_n\to Y_n$ is a weak equivalence in
$\topcat$ (i.e., a homotopy equivalence).
It is a \emph{cofibration} if for each $n\ge 0$, in the commutative square
\[
\xymatrix@R15pt{
\latch_n X \ar[d] \ar[r] & \ar[d] \latch_n Y \\
X_n   \ar[r] &  Y_n
}
\]
the induced map from the colimit of $X_n \leftarrow \latch_n X \to\latch_n Y$ to $Y_n$ is a cofibration
in $\topcat$. Here $\latch_nX\subset X_n$ is the $n$-the latching space,
\[  \latch_nX = \colimsub{[n]\twoheadrightarrow [p]} X_p  \]
where the colimit is taken over the category (poset) whose objects are surjective monotone maps $[n]\to [p]$,
for some $p<n$, and the morphisms are commutative triangles. It can also be regarded as a closed subspace of $X_n$,
the union of the images of all the simplicial operators $X_p\to X_n$ corresponding to surjective, non-bijective
monotone maps $[n]\to [p]$.

In particular, a simplicial space $Y$ is cofibrant if and only if the inclusions $\latch_nY\to Y_n$
are always (closed) cofibrations.

In this setting, there is a functorial (Reedy) fibrant replacement functor $\varphi$.
The formula for this is
\[ (\varphi X)_n := \holimsub{[m]\to [n] \textup{ in }\Delta_\inj} X_m  \]
where $[m]\to [n]$ runs over the comma category $(\Delta_\inj\!\downarrow\![n])$;
in particular $[m]\to [n]$ stands for an \emph{injective} monotone map.
The simplicial operators are defined as follows. Let $\varepsilon_n$ be the forgetful functor
from the comma category $(\Delta_\inj\!\downarrow\![n])$ to $\Delta$.
A morphism $f\co [m]\to [n]$ in $\Delta$ induces a functor
\[  f_!\co (\Delta_\inj\!\downarrow\![m])\to (\Delta_\inj\!\downarrow\![n])  \]
which is the composition of $f\circ \co (\Delta_\inj\!\downarrow\![m])\to (\Delta\!\downarrow\![n])$
with the left adjoint of the inclusion $(\Delta_\inj\!\downarrow\![n])\to (\Delta\!\downarrow\![n])$.
There is a preferred and obvious natural transformation
$\varepsilon_{m} \to \varepsilon_{n} f_!$
which induces contravariantly $X \varepsilon_{n} f_! \to X \varepsilon_{m}$ and then
\[ \rule{9mm}{0mm} (\varphi X)_{n}= \holim (X\varepsilon_{n}) \to \holim (X\varepsilon_{n}f_!) \to \holim (X\varepsilon_{m}) =(\varphi X)_{m} \,. \]

\begin{lem} \label{lem-manyreedycof} The simplicial spaces $\config(M)$, $\config(\pi)$, $\fin(M)$ and $\fin(\pi)$ are Reedy cofibrant. \qed
\end{lem}

\begin{lem} \label{lem-sermon}
Let $f\co X\to Y$ be a map of simplicial spaces in which $f$
is degreewise a covering space. Then the commutative square
\[
\xymatrix@R=20pt{
X \ar[d]_-f \ar[r] & \varphi X \ar[d]^-{\varphi f} \\
Y \ar[r] & \varphi Y
}
\]
is a strict pullback square, and both vertical arrows are degreewise covering spaces.
\end{lem}
\proof This is not about subtle properties
like being locally $1$-connected. By \emph{cove\-ring space} we mean a locally trivial map with discrete fibers.

Fix $m\ge 0$ and let $\iota_{Y,m}\co Y_m\to (\varphi Y)_m$ be the inclusion. It is a homotopy
equivalence and has a preferred left inverse $\rho_{Y,m}\co (\varphi Y)_m\to Y_m$.
Using similar notation for $X$, we can say that $\rho_{Y,m}$ and $\rho_{X,m}$ are the
horizontal arrows in a strict pullback square of spaces (in which the vertical arrows
are again $\varphi f$ and $f$). \qed

\medskip
For the remainder of this section we return to the setting of the introduction,
with manifolds $L$ and $M$ and covering spaces $\pi_L\co E_L\to L$ as well as $\pi_M\co E_M\to M$.

\begin{prop} \label{prop-wonderland} There is a homotopy pullback square, and in fact a strict pullback square,
of simplicial sets
\[
\xymatrix@R=20pt{
{\map_{\epifin\to\fin}
\left(\downmap{\config(\pi_L)}{\config(L)},\downmap{\varphi\config(\pi_M)}{\varphi\config(M)}\right)}
\ar[d] \ar[r] & \mapcov(\pi_L,\pi_M) \ar[d]^-q \\
\map_\fin(\config(L),\varphi\config(M)) \ar[r] & \map(L,M)
}
\]
\end{prop}
\proof This
needs a bit of unraveling. The right-hand column was already explained in the introduction. Here we allow ourselves to view it as a
map of simplicial sets (the singular simplicial sets of the spaces which we called $\mapcov(\pi_L,\pi_M)$ and $\map(L,M)$ in the introduction).
The top left-hand term is meant to be the space (simplicial set) of pairs $(f,g)$ making a commutative square
\[
\xymatrix@C=40pt{
\config(\pi_L) \ar[d]_-{g} \ar[r]^-{\textup{forgetful}} & \ar[d]^-f  \config(L) \\
\varphi\config(\pi_M) \ar[r]^-{\textup{forgetful}}  & \varphi\config(M)
}
\]
over the forgetful functor $\epifin\to \fin$.

For the proof we begin with the observation that we may replace $\varphi\config(M)$ and $\varphi\config(\pi_M)$
by $\varphi\fin_M$ and $\varphi\fin(\pi_M)$, respectively. This uses proposition~\ref{prop-configcov}.
(The functor $\varphi$ respects both homotopy pullback squares and strict pullback squares, so we can use a variant
of the proposition which has prefixes $\varphi$ everywhere.) Now let $X(L)=\config(\pi_L)$,
$X(M)=\config(\pi_M)$ and
\[  Y(L)=\config(L)\times_{\fin} \epifin, \qquad  Y(M)=\config(M)\times_{\fin} \epifin \]
so that there are forgetful maps $X(L)\to Y(L)$ and $X(M)\to Y(M)$.
There is no difference between $\varphi Y(M)$ and $\varphi\config(M)\times_{\fin}\epifin$. It remains to show that
the square
\begin{equation} \label{eqn-wonderland}
\begin{split}
\xymatrix@R=25pt@C=40pt{
{\map_\epifin\left(\downmap{X(L)}{Y(L)},\downmap{\varphi X(M)}{\varphi Y(M)}\right)}  \ar[d] \ar[r] & \mapcov(\pi_L,\pi_M) \ar[d]^-q \\
\map_\epifin(Y(L),\varphi Y(M)) \ar[r] & \map(L,M) \\
}
\end{split}
\end{equation}
is a strict pullback square. We make a few observations related to unique path lifting in covering spaces.
\begin{itemize}
\item[(a)] For every $r\ge 0$, the commutative squares
\[
\xymatrix@C=40pt@R=16pt@M=6pt{ X(L)_r \ar[r]^-{d_1d_2\cdots d_r} \ar[d]  &  X(L)_0 \ar[d] \\
Y(L)_r \ar[r]^-{d_1d_2\cdots d_r} &  Y(L)_0
}
\xymatrix@C=40pt@R=16pt@M=6pt{ X(M)_r \ar[r]^-{d_1d_2\cdots d_r} \ar[d]  &  X(M)_0 \ar[d] \\
Y(M)_r \ar[r]^-{d_1d_2\cdots d_r} &  Y(M)_0
}
\]
are strict pullback squares. (The vertical arrows are the forgetful maps, and they are covering spaces.
Remember that $d_1d_2\cdots d_r$ is the ultimate target operator.)
\item[(b)] The second of the pullback squares in (a) remains a pullback square if $\varphi$ is inflicted
on all the terms.
\end{itemize}
The proofs of (a) and (b) are by inspection. These statements imply that the forgetful maps $X(L)\to Y(L)$,
$X(M)\to Y(M)$ and $\varphi X(M)\to \varphi Y(M)$ are \emph{right fibrations}, i.e., Grothendieck
constructions associated with contravariant functors (in a homotopical sense) from $Y(L)$, $Y(M)$ and $\varphi Y(M)$ respectively
to the category of spaces. In addition the fibers here are discrete spaces, a.k.a.~sets, so that we are talking about
contravariant functors to the category of sets.

Let $\epifinone$ be the full subcategory of $\epifin$ spanned by the objects of the form $\uli k\twoheadrightarrow \uli 1$.
(This is of course isomorphic to the small category which has objects $\uli k$ for $k>0$ and all \emph{injective} maps as morphisms.)
Write $X^\epifinone(L)$, $Y^\epifinone(L)$ etc.~for the preimages of $N\epifinone$ in $X(L)$, $Y(L)$ etc.~under the reference
maps to $N\fin$. We can divide the square~\eqref{eqn-wonderland} into two by inserting a middle column
\[
\xymatrix@R=20pt{
{\map_{\epifinone}  \left(\downmap{X^\epifinone(L)}{Y^\epifinone(L)},\downmap{\varphi X^\epifinone(M)}{\varphi Y^\epifinone(M)}\right)} \ar[d] \\
{\map_\epifinone(Y^\epifinone(L), \varphi Y^\epifinone(M))}
}
\]
Of the two sub-squares which this creates, the one on the left is again a strict pullback square and a
homotopy pullback square, by inspection. (Loosely speaking this expresses the fact that the contravariant
functors on $Y(L)$, $Y(M)$ and $\varphi Y(M)$ mentioned above take certain ``disjoint unions'' to products.)
It remains to show that the other square
\[
\xymatrix@R=20pt{
{\map_{\epifinone}\left(\downmap{X^\epifinone(L)}{ Y^\epifinone(L)},\downmap{\varphi X^\epifinone(M)}{\varphi Y^\epifinone(M)}\right)}  \ar[d] \ar[r] & \mapcov(\pi_L,\pi_M) \ar[d]^-q \\
{\map_\epifinone(Y^\epifinone(L),\varphi Y^\epifinone(M))} \ar[r] & \map(L,M)
}
\]
so created is also a strict pullback square. But this is just a matter of unraveling the definitions and
simplifying the notation. For example, $Y^\epifinone(L)$ is $N\epifinone \times \textsf{P}L$, where $\textsf{P}L$ is the
path category of $L$ (whose nerve is weakly equivalent to the constant simplicial space $L$).
\qed

\section{Breaking the shackles}\label{sec-break}
See \cite[sec.8]{BoavidaWeissLong} and also \cite[sec.2]{BoavidaWeissconftensor} for the notion of a \emph{conservative map} of simplicial spaces or Segal spaces, and conservatization procedures.

\begin{thm} \label{thm-covconser} The forgetful functor $\config(\pi_L)\to \config(E_L)$ is a conservatization map
over $\fin$. More precisely, the standard reference map $\config(E)\to \fin$ is conservative and
the commutative triangle
\[
\xymatrix@C=2pt@M=6pt@R=15pt{
\config(\pi_L) \ar[rr] \ar[dr] &&  \config(E_L) \ar[dl]^-{\textup{ref. map}} \\
& \fin
}
\]
is a conservatization, i.e., is derived initial among maps over $\fin$ from
$\config(\pi_L)$ to simplicial spaces over $\fin$ which are conservative as such.
\end{thm}
\emph{Remark.} Write $\pi$ for $\pi_L$ and $E$ for $E_L$.
The commutative triangle in this theorem is a contraction of a commutative square:
\[
\xymatrix@C=48pt@M=6pt{
\config(\pi) \ar[r]^-{\textup{forgetful}} \ar[d]^-{\textup{ref.map}} &  \config(E) \ar[d]^-{\textup{ref.map}} \\
\epifin \ar[r]^-{\textup{functor \emph{source}}} & \fin
}
\]
The intuitive meaning of the theorem is roughly this. If we make morphisms in $\config(\pi)$ weakly
invertible whenever they are taken to invertible morphisms in $\fin$ under the diagonal map in the square,
then we obtain $\config(E)$, up to degreewise weak equivalence.

\medskip
Our proof of theorem~\ref{thm-covconser} relies on proposition~\ref{prop-betasuppress} below, a quotation
from~\cite{BoavidaWeissconftensor}.
The proposition refers to $\Lambda$, an endofunctor of the category of simplicial spaces over $N\sC$
which is left adjoint in a derived sense to the inclusion of the full subcategory of
conservative simplicial spaces over $N\sC$. (The appropriate notion of weak equivalence here
is degreewise weak equivalence of simplicial spaces over $N\sC$.) The functor $\Lambda$ was defined in
\cite[\S8]{BoavidaWeissLong} and is also used prominently in~\cite{BoavidaWeissconftensor}.

Let $\sC$ be a (discrete) small category. Let
\[
\xymatrix@C=6pt@R=15pt{ A \ar[dr] \ar[rr]^v && A' \ar[dl] \\
 &   N\sC
}
\]
be a commutative triangle of simplicial spaces and simplicial maps.
Suppose that $A$ and $A'$ are Segal spaces, both fiberwise complete over $N\sC$.

\begin{prop} \label{prop-betasuppress} If $A'$ is conservative over $N\sC$ and $v$ induces weak equivalences
from $(\Lambda A)_0$ to $(\Lambda A')_0$ and from $(\Lambda(A\downarrow x))_0$ to
$(\Lambda(A'\downarrow v(x)))_0$ for every $x\in A_0$~, then $v$ is a conservatization map.
\end{prop}
\emph{Remark.} Note that $(\Lambda A')_0\simeq A'_0$ and $(\Lambda(A'\downarrow v(x)))_0\simeq
(A'\downarrow v(x))_0$ by zigzags of weak equivalences, since $A'$ and $(A'\downarrow v(x))$ are
conservative over $N\sC$.

\proof[Proof of theorem~\ref{thm-covconser}] (This resembles a certain proof in \cite{BoavidaWeissconftensor}.)
We continue to write $E$ for $E_L$ and $\pi$ for $\pi_L$.
The plan
is to apply proposition~\ref{prop-betasuppress} with
$\sC=\fin$ and the Segal spaces $A=\config(\pi)$, $A'=\config(E)$. It is clear that $A'=\config(E)$ is fiberwise
complete over $N\fin$. It is also easy to verify that $A=\config(\pi)$ is fiberwise
complete over $N\fin$. (Since we have already noted that $\config(\pi)$ is fiberwise complete over
$\epifin$, it is enough to observe that $\epifin$ is fiberwise complete over $\fin$.)

For every $k\ge 0$, the space $Q=\emb(\uli k\,,E)$ is a stratified space.
Elements $f,g\co \uli k\to E$ of $\emb(\uli k\,,E)$ belong to the same stratum if and only if
\[  (\pi f(x)=\pi f(y)) \Leftrightarrow (\pi g(x)=\pi g(y)) \]
holds for all $x,y\in \uli k$; that is, if $\pi f$ and $\pi g$ define the same incidence relation. The set of strata is in
canonical bijection with the set of equivalence relations on $\uli k$\,. (Another way to say it: the set of
strata is in canonical bijection with the set of objects $\uli k_1\to \uli k_0$ in $\epifin$ such that $k_1=k$.)

We apply \cite[Corollary 9.3]{Miller2} or \cite[Corollary A.10.4]{Lurie} to this
stratification. So let $\sE\sP_Q$ be the exit path category associated with $Q$. This is a
topological category whose object space is the topological disjoint union of the strata of $Q$;
each stratum has the subspace topology inherited from $Q$. Then the inclusion of $\sE\sP_Q$ in the full
path category $\sP_Q$ of $Q$ leads to a weak equivalence
\[  |N\sE\sP^\op_Q| \lra |N\sP^\op_Q|~\simeq~Q. \]
That can also be written in the form of a weak equivalence
\[  \hocolimsub{[r]} (N\sE\sP^\op_Q)_r ~\lra~\hocolimsub{[r]} (N \sP^\op_Q)_r~. \]
Comparison with the definition of $\Lambda$, or of the variant $\Lambda^\flat$ also
described in \cite[\S8]{BoavidaWeissLong},
shows that this last map is the map of spaces
\[  (\Lambda A)_0 \lra (\Lambda A')_0  \]
induced by the forgetful map $A\to A'$ over $N\sC=N\fin$, where $A=\config(\pi)$ and $A'=\config(E)$ as before.
It follows that the first of the two hypotheses in proposition~\ref{prop-betasuppress} is satisfied.
The second of the two hypotheses turns into a special case of the first if we use corollary~\ref{cor-configcov}.
\qed

\section{Lifting maps between configuration categories}\label{sec-lifting}
We come to the proof of theorem \ref{thm-main} in a more precise formulation, theorem~\ref{thm-solution} below.
Instead of writing
\[
\xymatrix{
q^*\rmap_{\fin}(\config(L),\config(M)) \ar@{..>}[r] & \rmap_{\fin}(\config(E_L),\config(E_M))
}
\]
as in theorem~\ref{thm-main}, we can use lemma~\ref{lem-manyreedycof} and write
\[
\xymatrix{
q^*\map_{\fin}(\config(L),\varphi\config(M)) \ar@{..>}[r] & \map_{\fin}(\config(E_L),\varphi\config(E_M)).
}
\]

\begin{thm} \label{thm-solution} The following factorization problem in the shape of a commutative diagram
has a natural solution:
\[
\xymatrix@C=30pt@R=15pt{
 q^*\emb^t(L,M) \ar[d] \ar[rr] && \emb^t(E_L,E_M) \ar[d] \\
q^*\map_{\fin}(\config(L),\varphi\config(M)) \ar@{..>}[r] &? & \ar[l]_-{\simeq} \map_{\fin}(\config(E_L),\varphi\config(E_M))
}
\]
\end{thm}

\proof We substitute
$\map_{\fin}(\config(\pi_L),\varphi\config(E_M))$ for the question mark.
The lower right horizontal arrow is induced by the conservatization map of theorem~\ref{thm-covconser},
i.e., by the forgetful map from $\config(\pi_L)$ to $\config(E_L)$.
The lower left horizontal arrow is the composition of three,
\[
\xymatrix@M=8pt{
q^*\map_{\fin}(\config(L),\varphi\config(M)) \ar[d]_-{\cong}^-{\textup{Prop.}~\ref{prop-wonderland}}  &  \map_{\fin}(\config(\pi_L),\varphi\config(E_M))  \\
{\RR\map_{\epifin\shortrightarrow\fin}\left(\downmap{\config(\pi_L)}{\config(L)}, \downmap{\config(\pi_M)}{\config(M)}\right)} \ar[r]^-{\textup{forget}}
& \map_{\epifin}(\config(\pi_L),\varphi\config(\pi_M))  \ar[u]^-{\tau\circ}
}
\]
Here $\tau$ is the forgetful map from $\varphi\config(\pi_M)$ to $\varphi\config(E_M)$. \qed


\section{Localized configurations} \label{sec-loc}
Assume that $L$ and $M$ are \emph{smooth} manifolds.
One of the main points of \cite{BoavidaWeissLong} is a commutative square
\begin{equation} \label{eqn-oldsquare}
\begin{aligned}
\xymatrix@R=16pt{
\emb(L,M) \ar[r] \ar[d] & \rmap_\fin(\config(L),\config(M)) \ar[d] \\
\imm(L,M) \ar[r] & \rmap_\fin(\config^\loc(L),\config^\loc(M))
}
\end{aligned}
\end{equation}
where $\config^\loc(-)$ is a localized form of the configuration category. (Objects
of $\config^\loc(L)$ are objects of $\config(L)$ equipped with a morphism to a
singleton configuration.) If $\dim(M)-\dim(L)\ge 3$, this square is
homotopy cartesian \cite[Thm.1.1, Thm.5.1]{BoavidaWeissLong}. There is a variant for topological manifolds $L$, $M$
and topological embeddings
and immersions
\begin{equation} \label{eqn-topoldsquare}
\begin{aligned}
\xymatrix@R=16pt{
\emb^t(L,M) \ar[r] \ar[d] & \rmap_\fin(\config(L),\config(M)) \ar[d] \\
\imm^t(L,M) \ar[r] & \rmap_\fin(\config^\loc(L),\config^\loc(M))
}
\end{aligned}
\end{equation}
which is constructed in the same way. It is not known whether~(\ref{eqn-topoldsquare}) is always homotopy
cartesian for $\dim(M)-\dim(L)\ge 3$. But if in addition $L$ and $M$ are smooth (again) and $\dim(M)\ge 5$,
then the good properties
of~(\ref{eqn-oldsquare}) combined with a result of \cite{Lashof} imply that
we obtain a homotopy cartesian square from~(\ref{eqn-topoldsquare}) by
deleting certain path components in the left-hand column, namely, those which
do not contain smooth embeddings, respectively smooth immersions.

\begin{rem} \label{rem-T1} In diagram~\eqref{eqn-oldsquare}, the lower left-hand term $\imm(L,M)$ can be
replaced by $T_1\emb(L,M)$, the first Taylor approximation of $\emb(-,M)$ evaluated on $L$.
This makes a substantial difference only if
$\dim(L)=\dim(M)$ and $L$ has some compact components, for otherwise the inclusion $\emb(L,M)\to \imm(L,M)$,
viewed as a natural transformation of functors in the variable $L$, induces a weak equivalence
\[ T_1\emb(L,M)\to T_1\imm(L,M)\simeq
\imm(L,M). \]
Similarly, in diagram~\eqref{eqn-topoldsquare}, the lower left-hand term $\imm^t(L,M)$ can be
replaced by $T_1\emb^t(L,M)$, the first Taylor approximation of $\emb^t(-,M)$ evaluated on $L$.

The justification for this, both in the smooth setting and in the topological setting,
is the universal property of $T_1$ together with the observation that
the lower right-hand term in the square is 1-polynomial as a contravariant functor in the variable $L$.
(There is still no claim that the revised square, with $T_1\emb(L,M)$ instead of $\imm(L,M)$ or
$T_1\emb^t(L,M)$ instead of $\imm^t(L,M)$, is homotopy cartesian in the codimension $<3$ cases.)
\end{rem}

\medskip
Square~(\ref{eqn-topoldsquare}) can be viewed as a square of spaces over $\map(L,M)$.
Therefore, by base change along the notorious $q$, we get
\begin{equation} \label{eqn-topoldsquareslice}
\begin{aligned}
\xymatrix{
q^*\emb^t(L,M) \ar[r] \ar[d] & q^*\rmap_\fin(\config(L),\config(M)) \ar[d] \\
q^*\imm^t(L,M) \ar[r] & q^*\rmap_\fin(\config^\loc(L),\config^\loc(M))
}
\end{aligned}
\end{equation}
Now a new project emerges: to arrange the commutative square~(\ref{eqn-topoldsquareslice}) and
the square analogous to~(\ref{eqn-topoldsquare}) with $E_L$ and $E_M$ (in place of $L$ and $M$ respectively) in a commutative
cube, or at least a weaker (derived) variant of a commutative cube. This cube should of course
extend the commutative square of theorem~\ref{thm-solution} as well.

\smallskip
What we need, evidently, is a candidate for a Segal space which merits to be denoted
$\config^\loc(\pi)$, dependent on a covering space $\pi\co E\to M$.

\begin{defn} \label{defn-configcovloc} In degree $r$, the simplicial space $\config^\loc(\pi)$
is the preimage of the object $(\uli 1\to \uli 1)$ of $\epifin$ under
\[
\xymatrix@C=55pt{
(\config(\pi))_{r+1} \ar[r]^-{\textup{ultimate target}} & (\config(\pi))_0 \ar[r]^-{\textup{ref. functor}} & \ob(\epifin).
}
\]
where the ultimate target operator is $d_1d_2\cdots d_r d_{r+1}$.
The simplicial operators are defined in such a way that $d_i$ and $s_i$ in $\config^\loc(\pi)$ are the
appropriate restrictions of $d_{i+1}$ and $s_{i+1}$ in $\config(\pi)$.
\end{defn}
Intuitively, $\config^\loc(\pi)$ is the \emph{over category} associated
with $\config(\pi)$ and the space of objects taken to $(\uli 1\to \uli 1)$ by the reference functor.
Beware that the object
$(\uli 1\to \uli 1)$ of $\epifin$ is not a terminal object. Instead, an object $(f\co \uli k\to \uli \ell)$
admits a morphism to $(\uli 1\to \uli 1)$ in $\epifin$ if and only if $k=\ell$ and $f$ is the identity map.
In that case, the morphism to $(\uli 1\to \uli 1)$ is unique. An object of $\config(\pi)$ taken to $(\uli 1\to \uli 1)$
by the reference functor to $\epifin$ is nothing but a point $x\in E$ (and its value $\pi(x)\in M$, which is redundant).
Keeping all that in mind, we see that an object in $\config^\loc(\pi)$ consists of the following data:
\begin{itemize}
\item[(i)] a point $x\in E$;
\item[(ii)] a morphism $h$ in $\config(M)$ from some object, say $g\co \uli k\to M$ (injective map), to the singleton
configuration $\uli 1\mapsto \pi(x)\in M$,
\end{itemize}
and nothing more. Fortunately, by the unique path lifting property of $\pi$ there is a unique pair
consisting of an object $g^\lambda\co \uli k\to E$ in $\config(E)$ and a morphism $h^\lambda$ in $\config(E)$ from there
to the singleton configuration $\uli 1\mapsto x\in E$ such that $\pi\circ g^\lambda=g$ and $\pi\circ h^\lambda=h$.

\begin{prop} The forgetful map $\config(\pi)\to\config(E)$ induces a map from
$\config^\loc(\pi)$ to $\config^\loc(E)$ which is a degreewise weak equivalence of simplicial spaces over $\fin$.
\end{prop}
\proof The map as such is clear, and it is an inclusion map (allowing for the deletion of redundant data). Useful
observation: we are looking at a simplicial map over $M$. Moreover the projections from
$\config^\loc(\pi)$ to $M$ and from $\config^\loc(E)$ to $M$ (via $E$) are degreewise fibrations. For fixed $z\in M$,
the map of fibers over $z$ takes the form
\begin{equation} \label{eqn-locdecomp}
 \coprod_{x\in \pi^{-1}(z)} (\config(M)\downarrow z) \lra \coprod_{x\in \pi^{-1}(z)} (\config(E)\downarrow x).
\end{equation}
Here $z$, resp.~$x$, is shorthand for a one-point configuration $\uli 1\to M$ with image $\{z\}$,
resp.~$\uli 1\to E$ with image $\{x\}$. On the summand with label $x$, and in degree $r$, the map is given by lifting
all strings of $r+1$ composable morphisms in $\config(M)$ with ultimate target $z$ to strings of $r+1$ composable morphisms
in $\config(E)$ with ultimate target $x$.
This uses the unique path lifting property for $\pi\co E\to M$. In particular, the summand with label $x$ in the
left-hand side is mapped to the summand with the same label $x$ in the right-hand side.
Now it suffices to show that~\eqref{eqn-locdecomp} is a weak equivalence.
But this is obvious from \cite[Lem.4.1]{BoavidaWeissLong}. \qed

\medskip
Therefore definition~\ref{defn-configcovloc} leads to a refinement of theorem~\ref{thm-solution}.

\begin{thm} \label{thm-solutionloc} The following factorization problem in the shape of a commutative diagram
has a natural solution:
\[
\xymatrix@C=20pt@R=14pt{
 q^*\emb(L,M) \ar[d] \ar[rr] && \emb(E_L,E_M) \ar[d] \\
\ar[d] q^*\map_{\fin}(\config(L),\varphi\config(M)) \ar@{..>}[r] &? \ar@{..>}[d] & \ar[l]_-{\simeq} \map_{\fin}(\config(E_L),\varphi\config(E_M)) \ar[d] \\
q^*\map_{\fin}(\config^\loc(L),\varphi\config^\loc(M)) \ar@{..>}[r] &?? & \ar[l]_-{\simeq} \map_{\fin}(\config^\loc(E_L),\varphi\config^\loc(E_M))
}
\]
\end{thm}

\proof For the single question mark we substitute
$\map_{\fin}(\config(\pi_L),\varphi\config(E_M))$ as in the proof of theorem~\ref{thm-solution}.
For the double question mark we substitute
\[   \map_{\fin}(\config^\loc(\pi_L),\varphi\config^\loc(E_M)). \]
The broken arrow in the lower row is induced by the broken arrow in the middle row, which was
specified in the proof of theorem~\ref{thm-solution}. \qed

\medskip
The commutative diagram of theorem~\ref{thm-solutionloc} can be completed to a diagram in the shape of a
cube in a mechanical way. Because it is so mechanical, only some indications will be given. Deleting the middle row gives a diagram
\[
\xymatrix@C=20pt@R=20pt{
 q^*\emb(L,M) \ar[d] \ar[rr] && \emb(E_L,E_M) \ar[d] \\
q^*\map_{\fin}(\config^\loc(L),\varphi\config^\loc(M)) \ar@{..>}[r] &?? & \ar[l]_-{\simeq} \map_{\fin}(\config^\loc(E_L),\varphi\config^\loc(E_M))
}
\]
By the universal property of $T_1$ a new middle row can be inserted:
\[
\xymatrix@C=20pt@R=20pt{
 q^*\emb(L,M) \ar[d] \ar[rr] && \emb(E_L,E_M) \ar[d] \\
\ar[d] q^*T_1\emb(L,M) \ar[rr] && T_1 \emb(E_L,E_M) \ar[d] \\
q^*\map_{\fin}(\config^\loc(L),\varphi\config^\loc(M)) \ar@{..>}[r] &?? & \ar[l]_-{\simeq} \map_{\fin}(\config^\loc(E_L),\varphi\config^\loc(E_M))
}
\]
The amalgamation of this last diagram with the diagram of theorem~\ref{thm-solutionloc} (along the top and
bottom rows) is the cube.

\section{First order coherence in towers} \label{sec-cohere}
For this section, which is of a bureaucratic nature, we deviate slightly from the setup of section~\ref{sec-intro}.
We keep the topological manifolds $L$ and $M$.
In addition we make the following assumptions and introduce associated notation.
\begin{itemize}
\item[(i)] $M$ is connected and equipped with a \emph{universal} covering $A_M\to M$. Let $\Gamma$ be the
covering translation group of that. We like to view $A_M\to M$ as a principal $\Gamma$-bundle.
\item[(ii)] $L$ is contained in $M$ as a locally flat submanifold. Then the pullback (restriction) of the principal $\Gamma$-bundle
$A_M\to M$ to $L$ is a principal $\Gamma$-bundle on $L$, denoted $A_L\to L$.
\item[(iii)] Subgroups $\Gamma_1$ and $\Gamma_2$ of $\Gamma$ are specified such that $\Gamma_1\supset \Gamma_2$.
We may also write $\Gamma_0$ for $\Gamma$. Write $M(j):=A_M/\Gamma_j$ and $L(j):=A_L/\Gamma_j$ for $j=0,1,2$, so that
$L(0)\cong L$ and $M(0)\cong M$.
This gives us covering spaces
\[ M(k)\to M(j), \qquad L(k)\to L(j) \]
for $k,j\in\{0,1,2\}$ and $k>j$.
(Although we write $A_M/\Gamma_j$ and the like, the actions here are left actions.)
\end{itemize}
Let $q_j\co \map^{\Gamma_j}(A_L,A_M) \to \map(L(j),M(j))$ be the map obtained by passing to quotients for the
action of $\Gamma_j$. This is still a covering space. By making a base change along $q_j$ we obtain a covering space
\[   q_j^*\emb^t(L(j),M(j))\lra \emb^t(L(j),M(j)). \]
Now we can make a (commutative) diagram which is analogous to the diagram in theorem~\ref{thm-main}, but rotated by 90 degrees
and more repetitive:
\begin{equation} \label{eqn-tower1}
\begin{aligned}
\xymatrix@R=20pt{
q_2^*\emb^t(L(2),M(2)) \ar[r] &   q_2^*\rmap_\fin(\config(L(2)),\config(M(2))) \\
q_1^*\emb^t(L(1),M(1)) \ar[u] \ar[r] & \ar@{..>}[u]  q_1^*\rmap_\fin(\config(L(1)),\config(M(1)) \\
q_0^*\emb^t(L(0),M(0))             \ar[u] \ar[r] & \ar@{..>}[u]  q_0^*\rmap_\fin(\config(L(0)),\config(M(0)))
}
\end{aligned}
\end{equation}
The broken arrows are constructed like the broken arrow in theorem~\ref{thm-main}, and they represent
zigzags. (Their construction is not the main issue here.) On the other hand we can try to forget
$\Gamma_1$ and $L(1)$, $M(1)$. Then we still have covering spaces $L(2)\to L(0)$ and $M(2)\to M(0)$
etc., and using these we obtain by the method of theorem~\ref{thm-main} a commutative diagram
\begin{equation} \label{eqn-tower2}
\begin{aligned}
\xymatrix@R=20pt{
q_{2}^*\emb^t(L(2),M(2))  \ar[r] &   q_{2}^*\rmap_\fin(\config(L(2)),\config(M(2))) \\
q_{0}^*\emb^t(L(0),M(0)) \ar[u] \ar[r] & \ar@{..>}[u]  q_{0}^*\rmap_\fin(\config(L(0)),\config(M(0)))
}
\end{aligned}
\end{equation}
without making use of~\eqref{eqn-tower1}. The issue is whether~\eqref{eqn-tower2}
agrees with the square obtained by deleting the middle row of~\eqref{eqn-tower1}.
This is what we would call \emph{coherence in the tower}~\eqref{eqn-tower1}. The argument for coherence is quite predictable.

\smallskip
Let $S$ be a nonempty subset of $\{0,1,2\}$. We can arrange the elements in their order, $s_0<\cdots<s_r$
where $r\le 2$. Let $\pi_{M,S}$ and $\pi_{L,S}$
denote the diagrams
\[
\begin{array}{l}
M(s_0)\leftarrow M(s_1)\leftarrow \cdots \leftarrow M(s_r), \\
L(s_0)\leftarrow L(s_1)\leftarrow \cdots \leftarrow L(s_r)
\end{array}
\]
respectively. Each arrow in these diagrams is a covering space. If $S=\{j\}$ is a singleton,
then $\pi_{M,S}$ reduces to $M(j)$, and we prefer to write $M(j)$ in such a case.

What we need now is a definition of $\config(\pi_{L,S})$ which generalizes
definition~\ref{defn-epifin1} and a corresponding theorem which generalizes theorem~\ref{thm-covconser}.
We have all that for $|S|=2$, and for $|S|=1$ we can define $\config(\pi_{L,S}):=\config(L(j))$ where $j$ is the unique element of $S$.

\begin{defn} For $S=\{0,1,2\}$ let $\config(\pi_{L,S})$ be the nerve of the following topological category. An object is a
commutative diagram
\begin{equation} \label{eqn-multiconfcov}
\begin{aligned}
\xymatrix@C=12pt@R=14pt{
\uli k_0 \ar[d] & \ar[d] \ar[l] \uli k_1 & \ar[d] \ar[l] \uli k_2  \\
L(0) & \ar[l] L(1) & \ar[l] L(2)
}
\end{aligned}
\end{equation}
where the vertical arrows are injective, the arrows in the lower row are prescribed as in (iii) above, and the
arrows in the upper row are (surjective and) selfic maps.
A morphism from $(e_i\co \uli k_i\to L(s_i))_{i=0,1,2}$ to $(f_i\co \uli\ell_i\to L(s_i))_{i=0,1,2}$ is a family
$(\gamma_i)_{i=0,1,2}$ where $\gamma_i$ is a morphism (element of degree 1) in $\config(L(i))$, particle model, with source $e_i$ and target $f_i$.
\emph{Condition}: $\gamma_i$ is the unique multi-path lift of $\gamma_{i-1}$ with initial datum $e_i$ for the
covering space $L(s_i)\to L(s_{i-1})$, for $i=1,2$.
(Therefore $\gamma_2$ and $\gamma_1$ are determined by $\gamma_0$\,, $(e_i)$ and $(f_i)$, but their existence cannot be
taken for granted if only $(e_i),(f_i)$ and $\gamma_0$ are specified.)
\end{defn}
The definition suggests that we need another category for book-keeping. This should generalize $\epifin$ and we call it
$\trifin$. An object is a diagram $\uli k_0 \leftarrow \uli k_1\leftarrow\uli k_2$ where each arrow is (surjective and) selfic.
A morphism is a commutative diagram of maps of finite sets
\[
\xymatrix@C=12pt@R=14pt{
\uli k_0 \ar[d] & \ar[d] \ar[l] \uli k_1 & \ar[d] \ar[l] \uli k_2  \\
\uli \ell_0 &  \ar[l] \uli \ell_1 &  \ar[l] \uli \ell_2
}
\]
where the rows are objects of $\trifin$ and each of the little squares qualifies as a morphism in $\epifin$
(from $\uli k_{i-1}\leftarrow \uli k_i$ to $\uli\ell_{i-1}\leftarrow\uli \ell_i$). Then it is clear that for $S=\{0,1,2\}$ we have a
forgetful map of Segal spaces from $\config(\pi_{L,S})$ to (the nerve of) $\trifin$ induced by the functor taking an object~\eqref{eqn-multiconfcov}
to its upper row.

\begin{thm} \label{thm-covconsergen} For $S=\{0,1,2\}$, the forgetful map from $\config(\pi_{L,S})$
to $\config(L(2))$ is a conservatization map over $\fin$.
\end{thm}
This generalizes theorem~\ref{thm-covconser}. We make $\config(\pi_{L,S})$ into a simplicial space over $\fin$ by means of the composition
$\config(\pi_{L,S}) \to \trifin \to \fin$,
where the second arrow is the forgetful functor taking
$\uli k_0 \leftarrow \uli k_1\leftarrow \uli k_2$ to $\uli k_2$.

\proof This is identical with the proof of~\ref{thm-covconser} except for some renaming (in that $L(2)$ replaces
$E$ and $L(0)$ replaces $L$) and the following small but interesting change. For every
$k\ge 0$, the space $Q=\emb(\uli k,L(2))$ is a stratified space. Two elements $f,g\co \uli k\to L(2)$ of $Q$ belong
to the same stratum if and only if for $j=0,1$, the compositions of $f$ and $g$ with the preferred map
from $L(2)$ to $L(j)$ determine the same incidence relation on $\uli k$. The set of
strata is in canonical bijection with the set of objects $\uli k_0 \leftarrow \uli k_1\leftarrow \uli k_2$
in $\trifin$ such that $k_2=k$. \qed

\smallskip
Let us see how theorem~\ref{thm-covconsergen} implies coherence in the tower~\eqref{eqn-tower1}.
According to the standard instructions that tower must be implemented as in
\begin{equation} \label{eqn-cohere1}
\begin{aligned}
\xymatrix@R=14pt{
q_2^*\emb^t(L(2),M(2)) \ar[r] &   q_2^*\rmap_\fin(\config(L(2)),\config(M(2))) \ar[d]^-\simeq \\
                                            &           q_2^*\rmap_{\epifin\to\fin}(\config(\pi_{L,\{1,2\}}),\config(M(2)) \\
q_1^*\emb^t(L(1),M(1)) \ar[uu] \ar[r] & \ar[u]  q_1^*\rmap_\fin(\config(L(1)),\config(M(1)) \ar[d]^-\simeq \\
                                       &               q_1^*\rmap_{\epifin\to\fin}(\config(\pi_{L,\{0,1\}}),\config(M(1)) \\
q_0^*\emb^t(L(0),M(0))             \ar[uu] \ar[r] & \ar[u]  q_0^*\rmap_\fin(\config(L(0)),\config(M(0)))
}
\end{aligned}
\end{equation}
On the other hand diagram~\eqref{eqn-tower1} must be implemented as in
\begin{equation} \label{eqn-cohere2}
\begin{aligned}
\xymatrix@R=14pt{
q_2^*\emb^t(L(2),M(2))  \ar[r] &   q_2^*\rmap_\fin(\config(L(2)),\config(M(2))) \ar[d]^-\simeq \\
                                            &           q_2^*\rmap_{\epifin\to\fin}(\config(\pi_{L,\{0,2\}},\config(M(2)) \\
q_0^*\emb^t(L(0),M(0))             \ar[uu] \ar[r] & \ar[u]  q_0^*\rmap_\fin(\config(L(0)),\config(M(0)))
}
\end{aligned}
\end{equation}
The answer to our prayers for coherence is therefore
\[ K:= q_2^*\rmap_{\trifin\to\fin}(\config(\pi_{L,\{0,1,2\}}),\config(M(2))). \]
More precisely, in order to make the connection between~\eqref{eqn-cohere1} and~\eqref{eqn-cohere2}
we make alterations as in the sequence of schematic commutative diagrams
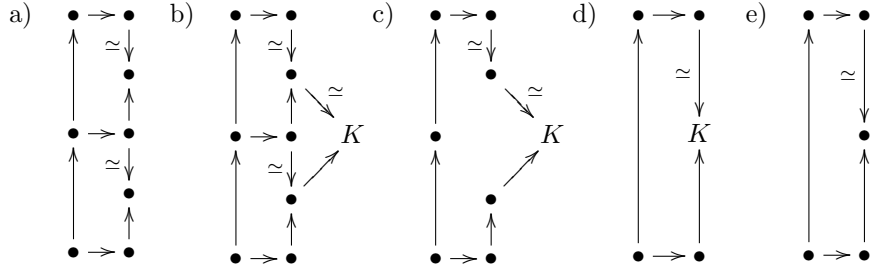
\[
\textup{a)}\quad
\xymatrix@R=12pt@C=10pt{
\bullet \ar[r] &  \bullet \ar[d]_-\simeq \\
                                            & \bullet  \\
\bullet  \ar[uu] \ar[r] & \ar[u]  \bullet \ar[d]_-\simeq \\
                                       &   \bullet  \\
\bullet             \ar[uu] \ar[r] & \ar[u]  \bullet
}  \quad \textup{b)}\quad
\xymatrix@R=12pt@C=10pt{
\bullet \ar[r] &  \bullet \ar[d]_-\simeq \\
                                            & \bullet \ar[dr]^-\simeq \\
\bullet  \ar[uu] \ar[r] & \ar[u]  \bullet \ar[d]_-\simeq & K  \\
                                       &   \bullet \ar[ur] \\
\bullet             \ar[uu] \ar[r] & \ar[u]  \bullet
}
\textup{c)}\quad
\xymatrix@R=12pt@C=10pt{
\bullet \ar[r] &  \bullet \ar[d]_-\simeq \\
                                            & \bullet \ar[dr]^-\simeq \\
\bullet  \ar[uu] &  & K  \\
                                       &   \bullet \ar[ur] \\
\bullet             \ar[uu] \ar[r] & \ar[u]  \bullet
}
\textup{d)}\quad
\xymatrix@R=14pt@C=10pt{
\bullet \ar[r] &  \bullet \ar[dd]_-\simeq \\
                                            &  \\
 &  K  \\
                                       &  \\
\bullet             \ar[uuuu] \ar[r] & \ar[uu]  \bullet
}
\quad\textup{e)}\quad
\xymatrix@R=14.5pt@C=10pt{
\bullet \ar[r] &  \bullet \ar[dd]_-\simeq \\
                                            &  \\
 &  \bullet  \\
                                       &  \\
\bullet             \ar[uuuu] \ar[r] & \ar[uu]  \bullet
}
\]
where a) stands for~\eqref{eqn-cohere1} and e) stands for~\eqref{eqn-cohere2}.
The move from a) to b)
uses the forgetful maps $\config(\pi_{L,\{0,1,2\}})\to \config(\pi_{L,\{1,2\}})$ and
$\config(\pi_{L,\{0,1,2\}})\to \config(\pi_{L,\{0,1\}})$. The first of these becomes a weak equivalence
after conservatization. The move from d) to e) uses the forgetful map
$\config(\pi_{L,\{0,1,2\}})\to \config(\pi_{L,\{0,2\}})$, which also becomes a weak equivalence after
conservatization.  --- In order to have something to refer to later,
we state this outcome as a proposition.
\begin{prop} Coherence holds in the tower~\eqref{eqn-tower1}. \qed
\end{prop}
It should be noted that theorem~\ref{thm-covconsergen} can also be used to establish first order coherence
in the more complicated setting of section~\ref{sec-loc}. This is for certain towers of the general shape
\[
\xymatrix@!@R=7pt@C=16pt{
& \bullet \ar@{->}[rr] \ar@{<-}'[d][dd]
& & \bullet \ar@{<..}[dd]
\\
\bullet \ar@{->}[rr]\ar@{->}[ur]\ar@{<-}[dd]
& & \bullet \ar@{->}[ur] \ar@{<..}[dd]
\\
& \bullet \ar@{->}'[r][rr]\ar@{<-}'[d][dd]
& & \bullet  \ar@{<..}[dd]
\\
\bullet \ar@{->}[rr]\ar@{->}[ur]\ar@{<-}[dd]
& & \bullet \ar@{->}[ur] \ar@{<..}[dd]
\\
& \bullet \ar@{->}'[r][rr]
& & \bullet
\\
\bullet \ar@{->}[rr]\ar@{->}[ur]
& & \bullet \ar@{->}[ur]
}
\]
with layers of the form
\[
\xymatrix@!C@C=-117pt@R=18pt@M=5pt{
& q_j^*T_1\emb^t(L(j),M(j))  \ar[rrr] &&&   q_j^*\rmap_\fin(\config^\loc(L(j)),\config^\loc(M(j))  \\
q_j^*\emb^t(L(j),M(j)) \ar[ur] \ar[rrr] &&& \ar[ur]  q_j^*\rmap_\fin(\config(L(j)),\config(M(j))
}
\]
where $j\in \{0,1,2\}$.


\begin{thebibliography}{10}
\bibitem{Andrade} R.~Andrade, \emph{From manifolds to invariants of $E_n$-algebras}, MIT thesis, 2010.
\bibitem{Boavida_Groth} P. Boavida de Brito, \emph{Segal objects and the Grothendieck construction},
 An alpine bouquet of algebraic topology, 19--44, Contemp.Math. vol. 708, Amer.Math.Soc., 2018.
\bibitem{BoavidaWeissHHA} P.~Boavida de Brito and M.S.~Weiss, \emph{Manifold calculus and homotopy sheaves}, Homology Homotopy Appl. 15 (2013), 361--383
\bibitem{BoavidaWeissLong} P.~Boavida de Brito and M.S.~Weiss, \emph{Spaces of smooth embeddings and configuration categories}, J. Topology 11 (2018), 65--143.
\bibitem{BoavidaWeissconftensor} P.~Boavida de Brito and M.S.~Weiss, \emph{The configuration category of a product},
 Proc.Amer.Math.Soc. 146 (2018), 4497--4512.
\bibitem{Cole2006} M.~Cole, \emph{Mixing model structures}, Topology and its Applications 153 (2006), 1016--1032.
\bibitem{DwyKa} W.~G.~Dwyer and D.~Kan, \emph{A classification theorem for diagrams of simplicial sets}, Topology 23 (1984), 139--155.
\bibitem{GoSche} P.~Goerss and K.~Schemmerhorn, \emph{Model categories and simplicial methods}, Interactions between
homotopy theory and algebra, 3--49, Contemp.Math.436, Amer.Math.Soc., 2007.
\bibitem{Kister} J.M.~Kister, \emph{Microbundles are fiber bundles}, Ann.of.Math. 80 (1964), 190--199.
\bibitem{Lashof} R.~Lashof, \emph{Embedding spaces}. Illinois J.Math. 20 (1976), no. 1, 144--154.
\bibitem{Lurie} J.~Lurie, \emph{Derived Algebraic Geometry VI: $E_k$ Algebras}, arXiv:0911.0018
\bibitem{Miller2} D.A.~Miller, \emph{Strongly stratified homotopy theory}, Tra.Amer.Math.Soc. 365 (2013), 4933--4962.
\bibitem{WeissEmb1} M.~Weiss, \emph{Embeddings from the point of view of immersion theory, I}, Geom.Top. 3 (1999), 67--101.
\end{thebibliography}
\end{document}